\documentclass[12pt]{amsart}

 \font\tenmsb=msbm10 scaled\magstep
1
\newfam\msbfam
   \textfont\msbfam=\tenmsb

\usepackage{amsmath}
\usepackage{amssymb}
\usepackage{mathscinet}
\usepackage{amsfonts,amsmath}
\usepackage{epsfig}
\usepackage{amsthm}
\usepackage{latexsym}
\usepackage[latin1]{inputenc}
\usepackage{verbatim}
\usepackage{pb-diagram}
\usepackage[small]{caption} 
\usepackage[labelformat=empty]{subfig}
\usepackage{everysel, keyval, ragged2e}
\usepackage{wrapfig}

\newtheorem{theorem}{Theorem}[section]

\theoremstyle{definition}
\newtheorem{definition}[theorem]{Definition}
\newtheorem{example}[theorem]{Example}

\theoremstyle{remark}

\numberwithin{equation}{section}

\newcommand{\bra}[1]{\langle #1 |}    
\newcommand{\ket}[1]{| #1 \rangle}    

\newcommand{\1}{{\rm 1\hspace*{-0.4ex}          
\rule{0.1ex}{1.52ex}\hspace*{0.2ex}}}

\newcommand{\m}{\mathcal }

\newcommand{\tr}{\operatorname{Tr} }

\def\C{\mathbb C}    

\begin{document}

\title[Private Channels, Conditional Expectations, Trace Vectors]{Private Quantum Channels, Conditional Expectations, and Trace Vectors}
\author[A.\ Church]{Amber Church$^1$}
\author[D.~W.\ Kribs]{David W. Kribs$^{1,2}$}
\author[R.\ Pereira]{Rajesh Pereira$^1$}
\author[S.\ Plosker]{Sarah Plosker$^1$}
\address{$^1$Department of Mathematics \& Statistics, University of Guelph, Guelph, ON, Canada N1G 2W1}
\address{$^2$Institute for Quantum Computing, University of Waterloo, Waterloo, ON, Canada N2L 3G1}

\begin{abstract}
Private quantum channels are the quantum analogue of the classical one-time pad. Conditional expectations and trace vectors are notions that have been part of operator algebra theory for several decades. We show that the theory of conditional expectations and trace vectors is intimately related to that of private quantum channels. Specifically we give a new geometric characterization of single qubit private quantum channels that relies on trace vectors. We further show that trace vectors completely describe the private states for quantum channels that are themselves conditional expectations. We also discuss several examples.
\end{abstract}
\maketitle
\section{Introduction}

Private quantum channels are a basic tool in quantum cryptography \cite{AMTW}. Conditional expectations and trace vectors are notions that have played a role in the theory of operator algebras for more than half a century \cite{Ume,MvN}. In this paper we show that there is an intimate relationship between the two subjects. Specifically we give a new geometric characterization of single qubit private quantum channels via the Bloch sphere representation for qubit states that relies on trace vectors. We further show that trace vectors completely describe the private states for quantum channels that are themselves conditional expectations.

In the next two sections we introduce private quantum channels, conditional expectations, and trace vectors. We discuss basic properties and include simple examples of each. We then consider the single qubit case in detail, giving a trace vector characterization of private states for unital quantum channels. We finish with a complete characterization of the private states for channels that are themselves conditional expectations in terms of trace vectors.

\section{Private Quantum Channels}\label{secQCandCE}

We will use $\m H$ or $\m K$ to denote Hilbert spaces (which are assumed to be finite dimensional in this paper) and denote a $d$-dimensional Hilbert space by $\m H_d$. We denote the set of linear operators on $\m H$ by $\m L(\m H)$, and we use $\mathbb M_d$ to denote the algebra of $d \times d$ complex matrices, which when convenient will be regarded as the matrix representations of operators in $\m L(\m H_d)$ with respect to a given orthonormal basis for $\m H_d$. The identity element of an operator space $X$ will be denoted by $\1_X$, or simply by $\1$ if the space is implied by the context, and we will write $\1_d$ for the identity of $\mathbb M_d$.

We will use Dirac notation for vectors $\ket{\phi}$ and vector duals $\bra{\phi}$. Thus pure states are represented as $\ket{\phi}\bra{\phi}$. General quantum states are represented by density operators (nonnegative operators with trace equal to 1), and we will use notation such as $\rho$, $\sigma$ in that case. The adjoint of an operator $A$ will be written $A^\dagger$, and we will reserve the asterisk notation when discussing abstract C$^*$-algebras.

A quantum channel is a linear, completely positive, trace preserving map $\m E:\m L(\m H)\rightarrow \m L(\m K)$. (Channels are generally defined on the set of trace class operators, with their dual maps defined on the set of bounded operators, but the sets coincide in the finite dimensional case.) Every channel can be written as
\begin{equation}\label{eq:Kraus}
\m E(\rho)=\sum_{i=1}^nK_i \rho K_i^\dagger,
\end{equation}
for some operators $K_i : \m H \rightarrow \m K$ with $\sum_{i=1}^nK_i^\dagger K_i=\1$, for any density operator $\rho$. We call a representation of $\m E$ as in equation (\ref{eq:Kraus}) a \emph{Kraus decomposition} of $\m E$. A channel $\m E$ is called a \emph{random unitary channel} if it admits a decomposition
\begin{equation}\label{eq:RUdecomp}
\m E(\rho)=\sum_ip_iU_i\rho U_i^\dagger \quad\quad\forall \rho,
\end{equation}
where $\{p_i\}$ form a probability distribution and $U_i$ are unitary operators.

In quantum cryptography, a private quantum channel (PQC) is the quantum analogue to the classical one-time pad. The following definition gives the mathematical framework for the notion in quantum information.

\begin{definition}
Let $\m S\subseteq \m H$ be a set of pure states and let $\m E:\m L(\m H)\rightarrow \m L(\m H)$ be a channel. Let $\rho_0$ be a density operator acting on $\m H$. Then $[\m S, \m E, \rho_0]$ is a \emph{private quantum channel} (PQC) if for any state $\ket{\phi}\in\m S$, we have
\[
\m E(\ket{\phi}\bra{\phi})=\rho_0.
\]
\end{definition}

PQCs were first considered in \cite{AMTW}, where the authors consider a particular class of random unitary channels. The most basic example is the following.

\begin{definition}
A map $\m E:\mathfrak{B}(\m H_{2^n})\rightarrow \mathfrak{B}(\m H_{2^n})$ is called a \emph{depolarizing channel} if, for any density matrix $\rho\in \mathfrak{B}(\m H_{2^n})$, we have
\[
\m E(\rho)=\frac{p}{n}\1+(1-p)\rho,
\]
where $0 < p \leq 1$ is a probability. When $p=1$ the \emph{completely depolarizing channel} is obtained, which gives the simplest example of a PQC, where every pure state is a private state. We denote the completely depolarizing channel by $\m E_{\C}$.
\end{definition}

\section{Conditional Expectations and Trace Vectors}\label{secCEandTV}

We recall a basic definition from operator algebras. Suppose there exists an orthogonal direct sum decomposition of a Hilbert space as $\m H=\bigoplus_i(M_i\otimes N_i)\oplus K$. Let $\m A$ be an algebra of operators in $\m L(\m H)$ consisting of all operators that belong to the set $\m A=\bigoplus_i(\1_{M_i}\otimes\m L(N_i))\oplus 0_K$, where $0_K$ is the zero operator on $K$. We call $\m A$ a \emph{(concrete finite dimensional) C$^*$-algebra}. $\m A$ is unital if $\1_{\m H}\in \m A$; i.e., $\m K$ is the zero subspace. A $*$-subalgebra $\m B$ of $\m A$ is a subset that is also a C$^*$-algebra. See \cite{Dav} for basic C$^*$-algebra theory.

\begin{definition} Let $\m A$ be a C$^*$-algebra and let $\m B\subseteq \m A$ be a unital $^*$-subalgebra. We call a linear map $\m E_{\m B}:\m A\rightarrow \m B$ a
\emph{conditional expectation} of $\m A$ onto $\m B$ if
\begin{enumerate}
\item[(i)] $\m E_{\m B}(b)=b$ for all $b\in \m B$;
\item[(ii)]$\m E_{\m B}(b_1ab_2)=b_1\m E_{\m B}(a)b_2$ for all $b_1, b_2\in \m B$ and for all $a\in \m A$;
\item[(iii)] $a \in \m A, \, a\ge 0$ implies $\m E_{\m B}(a)\ge 0$.
\end{enumerate}
\end{definition}

Conditional expectations were first considered in \cite{Ume}. We are interested in conditional expectations from $\mathbb{M}_n$ onto a subalgebra that are also quantum channels and hence trace preserving. We will therefore restrict ourselves to trace preserving conditional expectations.
We will call such maps \emph{conditional expectation channels}.

More examples of conditional expectation channels will be discussed below, but for the reader familiar with quantum information, we note here that the $n$-qubit completely depolarizing channel $\m E_{\C}$ is the conditional expectation onto the trivial scalar algebra ${\mathbb C} \cdot\1_{2^n}$. One way to see how conditional expectations inevitably arise in the theory is through trace inner products.

\begin{definition} A linear functional $\tau:\mathcal{A} \to \mathbb{C}$ is a \emph{faithful trace} if
\begin{enumerate}
\item[(i)] $\tau(a_1 a_2)=\tau(a_2 a_1)$
\item[(ii)] $\tau(a^\dagger a)>0$ for all $a\in \mathcal{A}$ with $a\neq 0$.
\end{enumerate}
\end{definition}

Given a faithful trace $\tau$ on $\mathcal{A}$ we can define an inner product $\langle a_1, a_2 \rangle=\tau(a_1^\dagger a_2)$.
We note that if $\mathcal{A}$ has a faithful trace $\tau$, the orthogonal projection onto $\mathcal{B}$ with respect to this inner product is the unique $\tau$-preserving conditional expectation from $\mathcal{A}$ to $\mathcal{B}$.  The essential structure of this argument can be found in \cite{Ume}. The most well-known example is the so-called Hilbert-Schmidt inner product $\langle A, B \rangle = \tr (A^\dagger B) $ on $\m M_n$.

\subsection{Trace Vectors}
We now consider trace vectors, a notion that initially arose in work of Murray and von Neumann \cite{MvN}, and has more recently been studied in the field of matrix theory.

\begin{definition}
Let $\m A$ be a $*$-subalgebra of $\m L(\m H_n)$. A vector $\ket{v}$ is a
\emph{trace vector} of $\m A$ if
\[
\bra{v}a\ket{v}=\frac1n\tr a \quad \forall a\in \m A.
\]
More generally, given a density operator $\rho_0$, we say $\ket{v}$ is a \emph{trace vector with respect to $\rho_0$} of $\m A$ if
\begin{equation}\label{eq:gentv}
\bra{v}a\ket{v}=\tr(\rho_0 a) \quad \forall a\in \m A.
\end{equation}
Thus by ``trace vector'', we really mean ``trace vector with respect to $\frac1n\1_n$''.
\end{definition}

By letting $a=\1$ in the definition of a trace vector, we find $\langle{v}|v\rangle=1$; a trace vector has unit length. It is easy to build a trace vector from other trace vectors in order to create a more general class of examples. Indeed, if $\ket{v_i}$ is a trace vector of the algebra $\m A_i=(\1_{M_i}\otimes\m L(N_i))\oplus 0_K$ for $i\in \{1,\dots, q\}$, then $\ket{v}=\bigoplus_{i=1}^q \ket{v_i}$ is a trace vector of the algebra $\m A=\bigoplus_{i=1}^q \m A_i$. In this way, trace vectors behave predictably. This also allows us to consider each summand separately, as we will do later.

\begin{example}\label{maxentangle}
As a fundamental example for quantum information, consider a maximally entangled state $\ket{\varphi_{e}}\in \m H_m\otimes \m H_n$. That is, a state
$\ket{\varphi_{e}}=\frac1{\sqrt{d}}\sum_{i=1}^d \ket{e_i}\otimes \ket{f_i}$, where $\{\ket{e_i}\}$ and $\{\ket{f_i}\}$ form orthonormal sets in $\m H_m$ and $\m H_n$ respectively, and $d=\min\{m,n\}$. If $m \geq n$, then one can check via direct calculation that $\ket{\varphi_{e}}$ is a trace vector for the algebra $\1_m \otimes \m L(\m H_n)$. And if $m = n$ an analogous calculation works for $\m L(\m H_m) \otimes \1_n$.
\end{example}

The general case is clarified by the following theorem of the third author. We recall that a vector $\ket{v}$ is a
\emph{separating vector} of an algebra $\m A$ if $a\ket{v}=0$ for some $a\in \m A$ implies $a=0$.

\begin{theorem}\label{rthm1}
If $\m A$ is a unital $*$-subalgebra of $\mathbb{M}_n$, then the following conditions are equivalent:
\begin{enumerate}
\item  $\m A$ is unitarily equivalent to
$
\oplus_{i=1}^q \left(\1_{m_i}\otimes \mathbb{M}_{n_i}\right),
$
where $m_i\ge n_i$ for all $i$ and $\sum_{i=1}^q m_in_i=n$.

\item $\m A$ has a separating vector.
\item $\m A$ has a trace vector.
\item There exists a set of trace vectors of $\m A$ that form an orthonormal basis of $\C^n$.
\end{enumerate}
\end{theorem}

This result is proved in \cite{Per03}. A related infinite dimensional open problem goes all the way back to von Neumann \cite{Ge}. Consider two simple cases. It is clear that $\mathbb{M}_n$ has no trace vectors -- from both the theorem and the definition of trace vectors. On the other hand, let $\Delta_2$ be the algebra of $2 \times 2$ diagonal matrices with respect to a basis $\{ \ket{0}, \ket{1}  \}$. One can readily check that the trace vectors for $\Delta_2$ are (up to complex phase) all vectors of the form $\ket{\psi} = \frac{1}{\sqrt{2}} ( \ket{0} + e^{i\theta} \ket{1} )$, for $0 \leq \theta < 2 \pi$; in other words, the set of all states that lie on the equator in the Bloch sphere representation for qubits (this point is further elucidated in the next section).

\section{Private Quantum Channels on the Bloch Sphere}
In this section we give a geometric characterization of single qubit unital PQCs in terms of the Bloch sphere representation \cite{NC} for single qubit states. We also show how the private states for such PQCs are determined by trace vectors. An alternative description was discussed in \cite{BZ}, where the entropy of sets of such private states was considered.

Every unital quantum channel is a random unitary channel in the single qubit case \cite{LS}. Thus, our private quantum channel $[\m S, \m E, \rho_0]$ in this case is given by a random unitary channel $\m E: \mathbb{M}_2\rightarrow \mathbb{M}_2$, a set of pure states $\m S$, and an output density matrix $\rho_0$. More precisely, we have $\m E(\ket{v}\bra{v})=\rho_0$ for all $\ket{v}\in \m S$. We would like to allow for the possibility of orthonormal vectors in $\m S$. As the channel is unital this can only occur if $\rho_0=\frac1n\1$, and hence we shall focus on this case here.

Using the Bloch sphere representation, we can associate
to any density matrix $\rho\in \mathbb{M}_2$ a Bloch vector $\vec{r}\in\mathbb{R}^3$ satisfying $\|\vec{r}\|\leq{1}$, where
\begin{equation}\label{eq:paulidecompofrho}
\rho=\frac{\1+\vec{r}\cdot    \vec{\sigma} }2.
\end{equation}
We use $\vec{\sigma}$ to denote the \emph{Pauli vector}, that is, $\vec{\sigma}=(\sigma_x,\sigma_y, \sigma_z)^T$.
Note that the set $\{\1, \sigma_x, \sigma_y, \sigma_z\}$ forms a basis for the \emph{real} vector space of Hermitian matrices in $\mathbb{M}_2$. We recall that a state is pure if and only if $\|\vec{r}\|=1$ and that the maximally mixed state $\frac {{\rm 1\hspace*{-0.3ex}\rule{.075ex}{1.14ex}\hspace*{0.15ex}}}{n}$ has Bloch vector $\vec{r}=\vec{0}$.

As discussed in \cite{KR01}, every linear map $\Phi:\mathbb{M}_2\rightarrow\mathbb{M}_2$ can be represented in the basis $\{\1, \sigma_x, \sigma_y, \sigma_z\}$ by a $4\times4$ matrix $\mathbb{T}$, and $\Phi$ preserves the trace if and only if the first row of the matrix $\mathbb{T}$ satisfies $t_{1k}=\delta_{1k}$; i.e.,
\begin{equation}
	\mathbb{T}=\left( \begin{array}{cc} 1 & \mathbf{0} \\
	\vec{t} & T \end{array} \right)
\end{equation}
where $T$ is a $3\times3$ matrix, $\mathbf{0}$ is a row vector, and $\vec{t}$ is a column vector.
The transformation $\Phi$ maps the subspace of Hermitian matrices into itself \emph{iff} $\mathbb{T}$ is real; finally, the map $\Phi$ is unital \emph{iff} $\vec{t}=\vec{0}$.

Thus, every unital qubit channel $\m E$ can be represented as
\begin{equation}\label{eq:channeldecomp}
	\m E\left(\frac1{2}\left[\1+\vec{r}\cdot\vec{\sigma}\right]\right)=\frac1{2}\left[\1+(T\vec{r})\cdot\vec{\sigma}\right],
\end{equation}
where $T$ and $\vec{t}$ are real, and we recall any density matrix can be written as in equation (\ref{eq:paulidecompofrho}). Here, the submatrix $T$ represents a deformation of the Bloch sphere, while the vector $\vec{t}$ represents a translation. This affine mapping of the Bloch sphere into itself is also discussed in section 8.3.2 of \cite{NC}.

We are of course interested in cases where $\m S$ is nonempty. This is easily seen to occur precisely when $T$ in equation (\ref{eq:channeldecomp}) has non-trivial nullspace.

Thus we consider the cases in which $T$ has non-trivial nullspace; that is, the subspace of vectors $\vec{r}$ such that $T \vec{r} = 0$ is one, two, or three-dimensional.

Finally, we note that in the single qubit case the unital subalgebras of the algebra $\m A=\mathbb{M}_2$ can be easily classified. They are
$\mathbb{M}_2$, $\C\cdot\1_2$ (the two trivial cases), and, up to unitary conjugation, $\Delta_2$, the subalgebra of all diagonal matrices in $\mathbb{M}_2$. To be precise, this third case refers to the subalgebras $\m B$ of the form $U^\dagger\Delta_2 U$, where $U\in \m A$ is unitary.

\begin{theorem}\label{thm:N(T)}
Let $\m E:\mathbb{M}_2\rightarrow\mathbb{M}_2$ be a unital qubit channel, with $T$ the mapping induced by $\m E$ as in equation (\ref{eq:channeldecomp}). Then there are three possibilities for a private quantum channel $[\m S, \m E, \frac12\1]$ with $\m S$ nonempty:
\begin{enumerate}
\item If the nullspace of $T$ is 1-dimensional, then $\m S$ consists of a pair of orthonormal states.
\item If the nullspace of $T$ is 2-dimensional, then the set $\m S$ is the set of all trace vectors of the subalgebra $U^\dagger\Delta_2 U$ of $2\times 2$ diagonal matrices up to a unitary equivalence.
\item If the nullspace of $T$ is 3-dimensional, then $\m E$ is the completely depolarizing channel and $\m S$ is the set of all unit vectors. In other words, $\m S$ is the set of all trace vectors of $\C\cdot \1_2$.
\end{enumerate}
\end{theorem}

\begin{figure}[h!]
  \centering
      \includegraphics[width=4in]{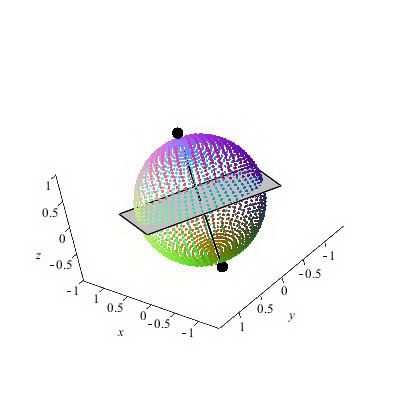}
  \vspace{-20pt}
  \caption{Case (1)}
\end{figure}

\begin{figure}[h!]
  \centering
      \includegraphics[width=4in]{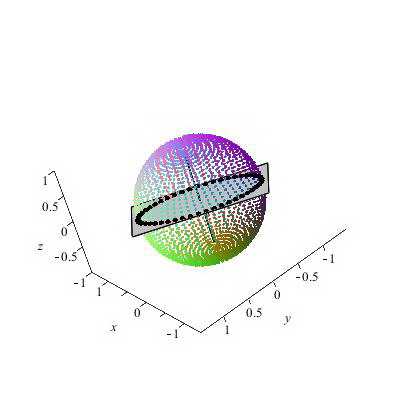}
  \vspace{-20pt}
  \caption{Case (2)}
\end{figure}


\proof
We shall write $\vec{r}_\phi$ for the Bloch sphere vector representation of a single qubit state $\ket{\phi}$. It is clear from equation~(\ref{eq:channeldecomp}) that $\m E(\ket{\phi}\bra{\phi}) = \frac{1}{2} \1$ if and only if $T  \vec{r}_\phi = 0$. Hence the relevant set that yields private states here is the intersection of the nullspace of $T$ and the surface of the Bloch sphere.

Case (1): The nullspace of $T$ is 1-dimensional. In this case, the nullspace is a single line through the origin of the Bloch sphere and the range of $T$ is a plane through the origin. Obviously this line meets the surface of the Bloch sphere in two antipodal points. These two antipodal points correspond to a pair of orthonormal single qubit states. Figure~1 gives an example.

Case (2): The nullspace of $T$ is 2-dimensional. In this case, the nullspace is a plane through the origin of the Bloch sphere. This plane meets the surface of the sphere in a great circle. The pure states corresponding to the points on this circle are precisely the private states for the channel. See an illustration in Figure~2.

To see how these private states arise from the trace vector perspective, let us consider the action of the channel more directly. As the nullspace of $T$ is 2-dimensional, its range is a line through the origin. For simplicity we shall assume this line is the $z$-axis; other cases are unitarily equivalent to this case. Thus, the range of $T$ intersects the sphere in the north and south poles, corresponding to the pure states $\ket{0}\bra{0}$ and $\ket{1}\bra{1}$ respectively.  The action of $T$ here will be a possible rotation of the Bloch sphere followed by a projection of the sphere onto the $z$-axis, followed by a possible contraction. By unitary equivalence, we only need consider the case where there is no initial rotation of the Bloch sphere.
In terms of the Pauli matrices $\sigma_x,\sigma_y,\sigma_z$, this means the action of the channel is given by $\m E(\sigma_x)=0$, $\m E(\sigma_y)=0$ and $\m E(\sigma_z) = p \sigma_z$ for some $0 < p \leq 1$.

Now $\Delta_2$ is the algebra of all diagonal matrices with respect to the ordered basis $\{ \ket{0}, \ket{1} \}$; explicitly, $\Delta_2$ is the set of all operators of the form $a \ket{0}\bra{0} + b \ket{1}\bra{1}$ for arbitrary scalars $a,b$. Then the projection onto the $z$-axis is a conditional expectation onto the subalgebra $\Delta_2$; call it $\m E_{\Delta}$. Explicitly,
\[
\m E_{\Delta}\left(
 \begin{bmatrix}
  a & b\\
  c & d
 \end{bmatrix}\right)= \begin{bmatrix}
  a & 0\\
  0 & d
 \end{bmatrix},
\text{ for any matrix } \begin{bmatrix}
  a & b\\
  c & d
 \end{bmatrix}.
 \]
One can check directly that $\m E = p \m E_{\Delta} + (1-p)\m E_{\C}$, where $\m E_{\C}$ is the completely depolarizing channel.

As $\m E_{\C}$ adds no restrictions to the private states for $\m E$, it suffices to show that the trace vectors for $\Delta_2$ are precisely the pure states that lie on the equator of the Bloch sphere. But the equator states are precisely the states that satisfy $|\bra{\phi}\ket{0}| = \frac{1}{\sqrt{2}} = |\bra{\phi}\ket{1}|$. And it is easy to see that these are the states which do indeed satisfy the trace vector condition for the algebra $\Delta_2$.

\begin{figure}[h!]
  \centering
      \includegraphics[width=3.5in]{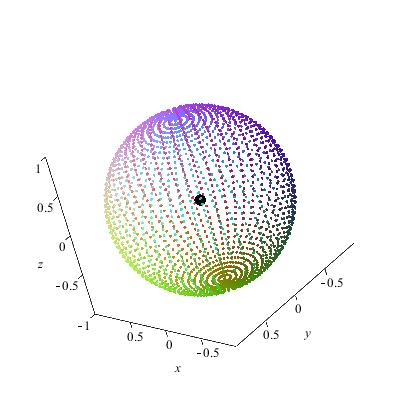}
  \vspace{-10pt}
  \caption{Case (3)}
\end{figure}

Case (3): The nullspace of $T$ is 3-dimensional, in other words $T$ is the zero operator. In this case $T$ maps the entire Bloch sphere to its origin, which corresponds to the maximally mixed state $\frac{1}{2} \1$, as shown in Figure~3. It is clear in this case that $\m E$ is the completely depolarizing channel $\m E_{\C}$. Moreover, the set $\m S$ has no restrictions; that is, $\m S$ is the set of all unit vectors. In other words, $\m S$ is the set of all trace vectors of $\C\cdot \1_2$.
\qed

\section{Private States for Conditional Expectation Channels}

The following result clarifies the general connection between conditional expectation channels, trace vectors and private states.

\begin{theorem}\label{thm:PQCiffTVs}
Let $\m E :\mathbb{M}_n \rightarrow\m A$ be a conditional expectation channel. Then $[\m S, \m E, \rho_0]$ is a private quantum channel if and only if $\m S$ is a set of trace vectors of $\m A$ with respect to $\rho_0 \in \m A$.
\end{theorem}

\proof
Let us first assume that $[\m S, \m E, \rho_0]$ is a PQC. Then $\m E(\ket{v}\bra{v}) = \rho_0$ for all $\ket{v} \in \m S$, and in particular note that $\rho_0$ belongs to $\m A$. Thus for all $\ket{v}\in \m S$ and for all $a \in \m A$, we have
\begin{eqnarray*}
\bra{v} a \ket{v}  &=&  \tr(\ket{v}\bra{v} a) \\
 &=&  \tr(\m E(\ket{v}\bra{v} a)) \\
 &=& \tr(\m E(\ket{v}\bra{v}) a)
= \tr(\rho_0 a),
\end{eqnarray*}
where the second and third identities follow from the trace preservation and conditional expectation properties of $\m E$ respectively. It follows that the states of $\m S$ are trace vectors of $\m A$ with respect to $\rho_0$.

For the converse, observe that when the vector states of $\m S$ are trace vectors of $\m A$ with respect to $\rho_0$, a similar calculation shows for all $\ket{v}\in\m S$ and for all $a\in \m A$ that
\begin{eqnarray*}
\tr(\rho_0 a)  &=&  \bra{v} a \ket{v} \\
 &=&  \tr(\ket{v}\bra{v} a) \\
 &=&  \tr(\m E(\ket{v}\bra{v} a))
 = \tr(\m E(\ket{v}\bra{v}) a).
\end{eqnarray*}
As $\rho_0$ belongs to $\m A$, it follows that $[\m S, \m E, \rho_0]$ forms a private quantum channel.
\qed

\begin{example}
Of course the three cases of Theorem~\ref{thm:N(T)} when applied to a unital single qubit conditional expectation channel $\m E_{\m A}:\mathbb{M}_2\rightarrow \m A$ are covered by this theorem. Indeed, applying Theorem \ref{thm:PQCiffTVs} to $\m E_{\m A}$ and letting $\rho_0=\frac12\1$ yields $[\m S, \m E_{\m A}, \frac12\1]$ is a PQC if and only if $\m A=U^\dagger\Delta_2U$ or $\m A=\C\cdot\1_2$ and $\m S$ is a set of trace vectors of $\m A$. Case 1 of  Theorem \ref{thm:N(T)} is an example of when $\m S$ is a proper subset of the set of all trace vectors of $U^\dagger\Delta_2U$, whereas Case 2 occurs when $\m S$ is the entire set. Case 3 occurs when $\m S$ is the set of all trace vectors of $\C\cdot\1_2$.
\end{example}

\begin{example}
Conditional expectations arise as the most basic non-trivial examples of private quantum communication using a private shared Cartesian frame \cite{BHS}. Let $\m H=(\C^2)^{\otimes N}$, and for simplicity suppose $N$ is even. Decompose the space as
\[
(\C^2)^{\otimes N}=\bigoplus_{j=0}^{N/2}\mathbb{H}_j\otimes \mathbb{K}_j,
\]
where the special unitary group $\operatorname{SU(2)}$ acts irreducibly on $\mathbb{H}_j$ and trivially on $\mathbb{K}_j$. As formulated in \cite{BHS}, if Alice and Bob share a reference frame to which Eve does not have access, and Alice prepares $N$ qubits in a state $\rho$ and sends them to Bob, Eve will see the resulting state simply as a mixture of all rotations $\Omega\in \operatorname{SU(2)}$. This situation can be summed up with the channel $\m E$, defined by
\[
\m E(\rho)=\sum_{j=0}^{N/2}(\m E_{\C j}\otimes id_{\mathbb{K}_j})(\Pi_j\rho\Pi_j),
\]
where $\m E_{\C j}$ is the completely depolarizing channel on $\mathbb{H}_j$ and $\Pi_j$ is the projection onto $\mathbb{H}_j$. One can see immediately that $\m E$ is in fact a conditional expectation channel that maps onto the algebra $\oplus_j (\1_{\mathbb{H}_j} \otimes \m L(\mathbb{K}_j))$. Thus, as noted in Theorem~\ref{thm:PQCiffTVs}, private states for $\m E$ can be found using trace vectors, which in this case can be constructed on the summands of the direct sum in a manner analogous to Example~\ref{maxentangle}.
\end{example}
\section{Outlook}

We see two main potential outcomes of the present work. Firstly, it is clear even just from the examples we have discussed here that there are numerous conditional expectation channels of relevance in quantum information, though they have not been viewed from this perspective before. It should be possible to use the conditional expectation and trace vector machinery to construct other new and useful examples of private quantum channels. Secondly, this work raises the intriguing possibility that a much more extensive theory of private quantum channels and private states could be developed. With few exceptions, the work on private channels appearing in the literature has focused primarily on specific instances and channels, rather than an overarching theory. We intend to continue these investigations elsewhere.

\vspace{0.1in}

{\noindent}{\it Acknowledgements.} We are grateful to Aron Pasieka for assistance with the Bloch sphere figures. D.W.K. was supported by Ontario Early Researcher Award 048142, NSERC Discovery Grant 400160 and NSERC Discovery Accelerator Supplement 400233. R.P. was supported by NSERC Discovery Grant 400096. S.P. was supported by an NSERC doctoral scholarship.

\bibliographystyle{amsalpha}

\end{document}